\documentstyle[amstex,amssymb,12pt]{article}
\newtheorem{Pa}{Paper}[section]
\newtheorem{Tm}[Pa]{{\bf Theorem}}
\newtheorem{La}[Pa]{{\bf Lemma}}
\newtheorem{Cy}[Pa]{{\bf Corollary}}
\newtheorem{Rk}[Pa]{{\bf Remark}}
\newtheorem{Pn}[Pa]{{\bf Proposition}}
\newtheorem{Pb}[Pa]{{\bf Problem}}
\newtheorem{Dn}[Pa]{{\bf Definition}}
\newtheorem{Ex}[Pa]{{\bf Example}}

\newcommand{\C}{{\mathbb C}}

\newcommand{\D}{{\mathbb D}}
\newcommand{\T}{{\mathbb T}}

\newcommand{\R}{{\mathbb R}}

\newcommand{\cH}{{\mathcal H}}

\newcommand{\cM}{{\mathcal M}}
\newcommand{\cN}{{\mathcal N}}

\newcommand{\cS}{{\mathcal S}}
\newcommand{\cZ}{{\mathcal Z}}
\newcommand{\cW}{{\mathcal W}}

\newcommand{\cQ}{{\Lambda}}

\newcommand{\sqminus}{{\rm{sq}_-}}
\newcommand{\eproof}{\hfill
    {\vbox{\hrule\hbox{\vrule height1.3ex\hskip0.8ex\vrule}\hrule}}\par}
\newcommand{\qed}{{\eproof}}

\begin{document}

\title{Functions with Pick matrices having bounded number of
negative eigenvalues}
\author{V. Bolotnikov, A. Kheifets and L. Rodman}
\date{}
\maketitle

\vskip 12pt

\begin{abstract}
A class is studied of complex valued functions defined on the unit disk
(with a possible
exception of a discrete set) with the property that all their Pick
matrices have not more than a prescribed number of negative eigenvalues.
Functions in this class, known to appear as pseudomultipliers
of the Hardy space, are
characterized in several other ways. It turns out that
a typical function in the class is meromorphic with a possible modification at
a finite
number of points, and total number of poles and of points of modification does
not exceed the prescribed number of negative eigenvalues.
Bounds are given for the number of points that generate Pick matrices that are
needed to achieve the requisite number of negative eigenvalues.
A result analogous to
Hindmarsh's theorem is proved for functions in the class.
\end{abstract}

\vskip 10pt
{\baselineskip=10pt
{\bf Key Words}: Schur functions, Pick matrices.}
\vskip 10pt
{\baselineskip=10pt
{\bf 2000 Mathematics Subject Classification}: 30E99, 15A63.}
\vspace*{.12in}

\baselineskip=15pt

\section{Introduction and main results}
\setcounter{equation}{0}

A complex valued function $S$ is called a {\em Schur function} if $S$ is
defined on the open unit disk $\D$, is analytic, and satisfies $|S(z)|\leq 1$
for every $z\in\D$. Schur functions and their generalizations have been
extensively studied in the literature. They admit various useful
characterizations; one such well-known characterization is
the following:
{\em A function $S$ defined on $\D$ is a Schur function
if and only
the kernel
\begin{equation}
K_S(z,w)=\frac{1-S(z)S(w)^*}{1-z{w}^*}
\label{1.2}
\end{equation}
is positive on} $\D$.
The positivity of (\ref{1.2}) on $\D$ means that
for every choice of positive integer $n$ and of distinct points
$z_1,\ldots,z_n\in\D$,
the {\em Pick matrix}
$$
P_n(S; \, z_1,\ldots,z_n)=
\left[\frac{1-S(z_i)S(z_j)^*}{1-z_i{z}^*_j}\right]_{i,j=1}^n
$$
is positive semidefinite: $P_n(S; \, z_1,\ldots,z_n)\geq 0.$
\smallskip

The following remarkable theorem due to Hindmarsh \cite{hind} (see also
\cite{D}) implies that if all Pick matrices of order $3$ are positive
semidefinite, then the function is necessarily analytic.
\begin{Tm}
Let $S$ be a function defined on an open set $U\subseteq \D$, and
assume that $P_3(S; \, z_1,z_2,z_3)$ is positive
semidefinite for every choice of $z_1, \, z_2, \, z_3\in U$. Then
$S$ is analytic on $U$, and $|S(z)|\leq 1$ for every $z\in U$.
\label{T:1.2}
\end{Tm}
A corollary of Theorem \ref{T:1.2} will be useful; we say that a set
$\cQ\subset \D$ is {\em discrete} (in $\D$) if $\cQ$ is at most countable,
with accumulation points (if any) only on the boundary of $\D$.
\begin{Cy}
Let $S$ be a function defined on $\D\setminus \cQ$, where $\cQ$ is a
discrete set. If $P_3(S; \, z_1,z_2,z_3)$ is positive
semidefinite for every choice of $z_1, \, z_2, \, z_3\in\D\setminus \cQ$,
then $S$ admits analytic continuation to a Schur function (defined on $\D$).
\label{T:1.2a}
\end{Cy}
In this paper we prove, among other things, an analogue of Hindmarsh's
theorem
for the class of functions whose Pick matrices have a bounded number
of negative eigenvalues.
These functions need not be analytic, or
meromorphic, on $\D$. It will be generally assumed that such functions
are defined on $\D\setminus \cQ$, where $\cQ$ is a discrete set, which
may depend on the function. More precisely:
\begin{Dn}
Given a nonnegative integer $\kappa$, the class $\cS_\kappa$ consists
of (complex valued) functions $f$ defined on ${\rm Dom}\,(f)=\D\setminus
\cQ$, where $\cQ=\cQ(f)$ is a discrete set, and such that all Pick
matrices
\begin{equation}
P_n(f; \, z_1,\ldots,z_n):=
\left[\frac{1-f(z_i)f(z_j)^*}{1-z_i{z}^*_j}\right]_{i,j=1}^n, \quad
z_1, \ldots ,z_n\in \D\setminus \cQ \ \ {\rm are} \ \ {\rm distinct},
\label{1.3a}
\end{equation}
have at most $\kappa$ negative eigenvalues, and at least one such Pick
matrix has exactly $\kappa$ negative eigenvalues (counted with
multiplicities).
\label{D:1.4}
\end{Dn}
The Pick matrices (\ref{1.3a}) are clearly Hermitian. Note also that in
Definition \ref{D:1.4} the points $z_1, \ldots ,z_n\in \D\setminus \cQ$
are assumed to be distinct; an
equivalent definition is obtained if we omit the requirement that the points
$z_1,\ldots,z_n$ are distinct.
\smallskip

Meromorphic functions in the class $\cS_{\kappa}$ have been
studied before in various contexts: approximation problems \cite{akh},
spectral theory of unitary operators in Pontryagin spaces \cite{kl},
Schur--Takagi problem \cite{AAK}, Nevanlinna - Pick problem
\cite{N}, \cite{gol}, \cite{BH}, model theory
\cite{ADRS}. Recently,
functions with jumps in
$\cS_{\kappa}$
appeared in the theory of almost multipliers (or pseudomultipliers)
\cite{AY}; this connection
will be discussed in more details later on.

Throughout the paper, ${\rm Dom}\,(f)$ stands for the domain of definition
of
$f$ and $Z(f)$ denotes the zero set for $f$:
$$
Z(f)=\{z\in{\rm Dom}\,(f): \; \; f(z)=0\}.
$$
The number of negative eigenvalues (counted with multiplicities) of
a Hermitian matrix $P$ will be denoted by $\sqminus P$.
For a function $f$ defined on $\D\setminus \cQ$, where $\cQ$
is a discrete set, we let ${\bf k}_n(f)$ to denote the maximal number
of negative eigenvalues (counted with multiplicities) of all Pick
matrices of order $n$:
\begin{equation}
{\bf k}_n(f):=\max_{z_1,\ldots,z_n\in\D\setminus \cQ}
\sqminus P_n(f; \, z_1,\ldots,z_n).
\label{1.4}
\end{equation}
We also put formally ${\bf k}_0(f)=0$. It is clear that the sequence
$\{{\bf k}_n(f)\}$ is non decreasing and
if $f\in\cS_\kappa$, then
\begin{equation}
{\bf k}_n(f)=\kappa
\label{1.5}
\end{equation}
for all sufficiently large integers $n$. Note that the class $\cS_0$
coincides with the Schur class; more precisely,
every function $f$ in $\cS_0$
admits an extension to a Schur
function. Here and elsewhere, we say that a function $g$ is an
{\em extension} of a function $f$ if ${\rm Dom}\,(g)\supseteq {\rm
Dom}\,(f)$
and $g(z)=f(z)$ for every $z\in {\rm Dom}\,(f)$.

\smallskip

The following result by Krein - Langer \cite{kl} gives a characterization
of
{\em meromorphic} functions in $\cS_\kappa$.
\begin{Tm}
Let $f$ be a meromorphic function on $\D$. Then $f\in\cS_\kappa$ if and
only if $f(z)$ can be represented as $f(z)=
{\displaystyle \frac{S(z)}{B(z)}},$ where $S$ is a
Schur function and $B$ is a Blaschke product of degree $\kappa$ such that
$S$ and $B$ have no common zeros.
\label{T:1.2b}
\end{Tm}
Recall that a {\em Blaschke product} (all Blaschke products in this paper are
assumed to be finite) is a rational function $B(z)$ that is
analytic on $\D$ and unimodular on the unit circle $\T$ : $|B(z)|=1$ for
$|z|=1$; the {\em degree} of $B(z)$ is the number of zeros (counted with
multiplicities) of $B(z)$ in $\D$.
See also \cite{gol}, \cite{DLS}, \cite{BR} for various proofs of matrix
and
operator--valued versions of Theorem \ref{T:1.2b}.
\smallskip

However, not all functions in $\cS_{\kappa}$ are meromorphic, as a
standard example shows:
\begin{Ex}
{\rm Let the function $f$ be defined on $\D$ as $f(z)=1$ if $z\neq 0$;
$f(0)=0$.
Then $P_n(f;z_1, \ldots ,z_n)=0$ if $z_j$'s are all nonzero; if one of
the points is zero (up to a permutation similarity we can assume that
$z_1=0$), then $P_n(f;z_1, \ldots ,z_n)$ is the matrix with ones in
the first column and in the first row and with zeroes elsewhere. This
matrix clearly is of rank two and has one negative eigenvalue.
Thus, $f\in \cS_{1}$. That $f$ has a jump discontinuity at $z=0$ is
essential; removing the discontinuity
would result in the constant function $\tilde{f}(z)=1$,
which does not belong to} $\cS_1$.
\label{E:2.1}
\end{Ex}
As it was shown in \cite{AY}, functions in $\cS_{\kappa}$
that are defined in $\D$ except for a finite set of singularities,
are exactly the $\kappa$-pseudomultipliers of the Hardy space $H_2$
(see \cite{AY} for the definitions of pseudomultipliers and their
singularities; roughly speaking, a pseudomultiplier maps a subspace
of finite codimension  in Hilbert space of functions into the Hilbert
space, by means of the corresponding multiplication operator). In fact,
similar results were obtained in \cite{AY} for more general reproducing
kernel Hilbert spaces of functions. Theorem 2.1 of \cite{AY} implies in
particular that every function in the class $\cS_{\kappa}$ defined on a
set of uniqueness of $H_2$ can be extended to a function in $\cS_{\kappa}$
defined everywhere in $\D$ except for a finite number of singularities.
Theorem 3.1 of \cite{AY} implies that the pseudomultipliers of $H_2$ are
meromorphic in $\D$ with a finite number of poles and jumps, some of which
may coincide. This result applies to other reproducing kernel Hilbert
spaces of analytic functions as well, as proved in \cite{AY}.
\smallskip

For the Hardy space $H_2$, the pseudomultipliers can be characterized in a
more concrete way, namely as standard functions defined below, see Theorem
\ref{T:1.4}(2). This characterization is based
largely on Krein - Langer theorem \ref{T:1.2b}.
\smallskip

In this paper we focus on a more detailed study of the class
$\cS_{\kappa}$. Our proofs depend on techniques used in
interpolation, such as matrix inequalities and Schur complements,
and allow us to obtain generalizations of Hindmarsh's theorem
(Theorem \ref{T:1.4}(3))
and precise
estimates of the size of Pick matrices needed
to attain the required number of
negative eigenvalues (Theorem \ref{T:3.1}).

\begin{Dn}
A function $f$ is said to be a {\em standard function}
if it admits the representation
\begin{equation}
f(z)=\left\{\begin{array}{cl} {\displaystyle\frac{S(z)}{B(z)}}
& \quad \mbox{if } \; z\not \in \cW\cup\cZ, \\
\gamma_j & \quad \mbox{if } \; z=z_j\in\cZ, \end{array}\right.
\label{1.6a}\end{equation}
for some complex numbers $\gamma_1 ,\ldots , \gamma_{\ell}$, where:
\begin{enumerate}
\item $\cZ=\{z_1, \ldots ,z_{\ell}\}$ and $\cW=\{w_1,\ldots ,w_p\}$ are
disjoint sets of distinct points in $\D$;
\item $B(z)$ is a Blaschke product of degree $q\geq 0$ and $S(z)$ is a
Schur function with the zero sets $Z(B)$ and $Z(S)$, respectively, such
that
$$
\cW\subseteq Z(B)\subseteq \cW\cup\cZ\quad\mbox{and}\quad
Z(B)\cap Z(S)=\emptyset;
$$
\item if $z_j\in \cZ\setminus Z(B)$, then $
{\displaystyle \frac{S(z_j)}{B(z_j)}}\neq \gamma_j.$
\end{enumerate}
\label{D:1.7}
\end{Dn}
For the standard function $f$ of the form (\ref{1.6a}), ${\rm Dom}\,(f)=
\D\setminus \cW$.
More informally, $\cZ$ is the set of points $z_0$ at which
$f$ is defined and $f(z_0)\neq \lim_{z \rightarrow z_0} f(z)$.
The case when $\lim_{z \rightarrow z_0} |f(z)|=\infty$ is not excluded here; in
this case $f(z_0)$ is defined nevertheless.
The set $\cW$
consists of those poles of $
{\displaystyle \frac{S(z)}{B(z)}}$ at which $f$ is not defined.
In reference to the properties (1)-(3) in Definition
\ref{D:1.7} we will say that $f(z)$ has $q$ poles (the zeros of $B(z)$),
where each pole is counted according to its multiplicity as a zero of $B(z)$,
and $\ell$ jumps $z_1, \ldots ,z_{\ell}$. Note that the poles and jumps
need not be disjoint.
\begin{Tm}
Let $f$ be defined on $\D\setminus \cQ$, where $\cQ$ is a discrete set.
Fix a nonnegative integer $\kappa$. Then the following statements are
equivalent:
\begin{enumerate}
\item $f$ belongs to $\cS_\kappa$.
\item $f$ admits an extension to a
standard function with $\ell$ jumps (for
some $\ell$, $0\leq \ell\leq \kappa$) and
$\kappa - \ell$ poles, and
the jumps of the standard function
are contained in $\D\setminus \cQ$.
\item \begin{equation} {\bf k}_n(f)={\bf k}_{n+3}(f)=\kappa
\quad  {\rm for} \quad {\rm  some} \quad  {\rm integer} \ \  n\geq 0.
\label{2.10}
\end{equation}
\end{enumerate}
\label{T:1.4}
\end{Tm}
Note that the extension to a standard function in the second statement is
unique. Note also that the equivalence $ {\bf 1 \Leftrightarrow 3}$ is a
generalization of Hindmarsh's theorem
to the class of functions whose Pick matrices have a bounded number of
negative eigenvalues,
and
coincides with Theorem \ref{T:1.2} if $\kappa=0$.
\smallskip

The third statement in
Theorem  \ref{T:1.4} raises the question of the minimal possible
value of $n$ for which (\ref{1.5}) holds. For convenience of future
reference we denote this minimal value by $N(f)$:
\begin{equation}
N(f)=\min_{n}\{n: \; {\bf k}_{n}(f)=\kappa\}, \qquad f\in\cS_{\kappa}.
\label{1.9}
\end{equation}
\begin{Tm}
Let $f$ be a standard function with $q$ poles and $\ell$
jumps. Then $f\in \cS_{\kappa}$, where $\kappa=q+\ell$, and
\begin{equation}
q+\ell\leq N(f)\leq q +2\ell.
\label{3.1}
\end{equation}
\label{T:3.1}
\end{Tm}
The left inequality in (\ref{3.1}) is self-evident (to have $\kappa$
negative eigenvalues, a Hermitian matrix must be of size at least
$\kappa\times
\kappa$), while the right inequality is
the essence of the theorem.
We shall show
that inequalities (\ref{3.1}) are exact; here we present a simple example
showing that these inequalities can be strict.
\begin{Ex}\label{E:last} {\rm Let
$$
f(z)=\left\{\begin{array}{cc} 1, & z\neq 0, \\
a\neq 1, & z=0.\end{array}\right.
$$
Then $f(z)$ is a standard function from $\cS_1$ with no
poles and one jump and Theorem \ref{T:3.1} guarantees that
$1\leq N(f)\leq 2$. More detailed analysis shows that if $|a|>1$, then
${\bf k}_n(f)=1$ for $n\ge 1$ and therefore $1=N(f)<2$. If $|a|\leq 1$,
then ${\bf k}_1(f)=0$, ${\bf k}_n(f)=1$
for $n\ge 2$ and therefore $1<N(f)=2$.}
\end{Ex}
The following addition to Theorem \ref{T:3.1} is useful:
\begin{Rk}
Let $f$ be as in
Theorem \ref{T:3.1}.
Let $w_1, \ldots ,w_k$ be the distinct poles of $f$ with multiplicities
$r_1, \ldots ,r_k$, respectively, and let $z_1,\ldots ,z_{\ell}$ be the
(distinct) jumps
of $f$. Then there exists an $\epsilon>0$ such that
every set
$Y=\{y_1,\ldots,y_m\}$ of
$m:=q+2\ell$ distinct points
satisfying the conditions 1 - 3 listed below
has the property that
$$
\sqminus P_m(f; y_1,\ldots ,y_m)=q+\ell.$$
\begin{enumerate}
\item $\ell$ points in the set $Y$ coincide with $z_1,\ldots,
z_{\ell}$;
\item
$\ell$ points in $Y$, different from $z_1,\ldots,
z_{\ell}$, are at a
distance less than $\epsilon$
from
the
set $\{z_1, \ldots, z_{\ell}\}$;
\item the remaining $q$ points in $Y$ are
distributed
over $k$ disjoint subsets $Y_{1}, \ldots, Y_k$ such that the set $Y_j$
consists
of $r_j$ points all at a positive distance less then $\epsilon$ from
$w_j$, for
$j=1,\ldots ,k$.
\end{enumerate} \label{R:last}
\end{Rk}
The proof of Remark \ref{R:last} is
obtained in the course of the proof of Theorem \ref{T:3.1}.
\smallskip

Still another characterization of the class $\cS_{\kappa}$ is given in the
following theorem:
\begin{Tm} Let $f$ be as in Theorem $\ref{T:1.4}$, and fix a nonnegative
integer $\kappa$. Then:
\begin{enumerate}
\item If
${\bf k}_{2\kappa}(f)={\bf k}_{2\kappa+3}(f)=\kappa$, then $f\in \cS_{\kappa}$.
\item If $f\in \cS_{\kappa}$, then ${\bf k}_{2\kappa}=\kappa.$
\end{enumerate}
\label{T:1.new}
\end{Tm}
Example \ref{E:last} (with $|a|\le 1$) shows that for a fixed
$\kappa$,
the subscript $2\kappa$ in Theorem \ref{T:1.new}
cannot be replaced by any smaller nonnegative integer.
\smallskip

Theorems \ref{T:1.4}, \ref{T:3.1}, and \ref{T:1.new} comprise the main results
of this paper.

\smallskip

Sections 2 through 4 contain the proofs of Theorems \ref{T:1.4} and
\ref{T:3.1}. In Section 5, a local result is proved to the effect that
for a function $f\in\cS_{\kappa}$ and any open set $\Omega$ in
${\rm Dom}\, (f)$, the requisite number of negative eigenvalues of Pick
matrices $P_n(f;z_1, \ldots, z_n)$ can be achieved when the points $z_j$
are restricted to belong to $\Omega$. Hindmarsh sets and their elementary
properties are introduced in Section 6, where we also state an open
problem.

\smallskip

If not stated explicitly otherwise, all functions are
assumed to be scalar (complex
valued). The superscript $^*$ stands for the conjugate of a complex
number or the conjugate transpose of a matrix.

\section{Theorem \ref{T:1.4}: part of the proof}
\setcounter{equation}{0}

In this section we prove implication ${\bf 3 \Rightarrow 1}$ of Theorem
\ref{T:1.4}. We start with some auxiliary lemmas.
\begin{La}
If an $m \times m$ matrix $X$ has rank $k<m$ and the algebraic
multiplicity of zero as an eigenvalue of $X$ is $m-k$, then
there exists a $k \times k$ nonsingular principal submatrix of $X$.
\label{L:4.1}
\end{La}
{\bf Proof:} The characteristic polynomial of $X$ has the form
$\lambda^m +a_{m-1}\lambda^{m-1}+\cdots + a_{m-k}\lambda^{m-k}$,
where $a_{m-k}\neq 0$. Since $\pm a_{m-k}$ is the sum of all
determinants of $k \times k$ principal submatrices of $X$, at least
one of those determinants must be nonzero.  \qed
\begin{La}
\label{L:new}
Let $X$ be a Hermitian $n\times n$ matrix. Then
\begin{enumerate}
\item $\sqminus Y\geq\sqminus X$ for all Hermitian matrices
$Y\in\C^{n\times n}$ in some small neighborhood of $X$.
\item $\sqminus P^*XP \leq \sqminus X$ for every $P\in\C^{n\times m}$.
In particular, if $X_0$ is any principal submatrix of $X$, then
$\sqminus X_0 \leq \sqminus X$.
\item If $\{X_m\}_{m=1}^\infty$ is a sequence of Hermitian matrices such that
$\sqminus X_m\leq k$ for $m=1,2, \ldots$, and $\lim_{m\rightarrow \infty}
X_m=X$, then also $\sqminus X\leq k$.
\end{enumerate}
\end{La}
{\bf Proof}: The first and third statements easily follow by the
continuity properties of eigenvalues of Hermitian matrices.
The second statement is a consequence from the interlacing properties of
eigenvalues of Hermitian matrices, see, e.g. \cite[Section 8.4]{LT}. \qed
\begin{La} \label{L:last1}
Let $f$ be a function defined on $\D\setminus \cQ$, where $\cQ$
is a discrete set.
\begin{enumerate}
\item If $g$ is defined on a set $(\D\setminus \cQ)\cup
\{w_1,\ldots,w_k\}$
(the points $w_1,\ldots,w_k\in\D$ are not assumed to be distinct or to be
disjoint with $\cQ$) and $g(z)=f(z)$ for every $z\in \D\setminus \cQ$,
then
\begin{equation}
{\bf k}_n(f)\leq {\bf k}_n(g) \leq
\max_{0\leq r\leq \min \{k,n\}} \{{\bf k}_{n-r}(f) +r\}
\leq
{\bf k}_n(f) +k, \qquad n=1,2,
\ldots.
\label{4.1a}
\end{equation}
\item If $g$ is defined on $\D\setminus \cQ$, and $g(z)=f(z)$ for every
$z\in \D\setminus \cQ$, except for $k$
distinct points $z_1, \ldots z_k\in
\D\setminus \cQ$, then
$$
{\bf k}_n(g) \leq  {\bf k}_n(f)+k,
\qquad n=1,2,
\ldots.
$$
\end{enumerate}
\end{La}
{\bf Proof:} We prove the first statement.
The first inequality in (\ref{4.1a}) follows from the definition
(\ref{1.4}) of
${\bf k}_n(f)$. Let $y_1, \ldots ,y_n$ be a set of
distinct points in $(\D\setminus \cQ)\cup
\{w_1,\ldots,w_k\}$, and assume that exactly $r$ points among
$y_1, \ldots ,y_n$ are in the set $\{w_1, \ldots w_k\}\setminus
(\D\setminus \cQ)$. For notational convenience, suppose that
$y_1, \ldots ,y_r \in \{w_1, \ldots w_k\}\setminus
(\D\setminus \cQ)$. We have
$$ P_n(g;y_1, \ldots ,y_n)=\left[\begin{array}{cc} ^* & ^* \\ ^* &
P_{n-r}(f;y_{r+1},\ldots ,y_n) \end{array}\right], $$
where $*$ denotes blocks of no immediate interest. By the
interlacing properties of eigenvalues of Hermitian matrices, we have
$$\sqminus P_n(g;y_1, \ldots ,y_n)\leq
\sqminus P_{n-r}(f;y_{r+1},\ldots ,y_n) +r \leq {\bf k}_{n-r}(f) +r, $$
and, since $0\leq r\leq \min\{n,k\}$, the second and third inequalities in
(\ref{4.1a}) follow.

For the second statement, using induction on $k$, we need to prove only the
case $k=1$. Let $y_1, \ldots ,y_n$ be a set of
distinct points in $\D\setminus \cQ$. Then
\begin{equation}\label{ag} P_n(g;y_1, \ldots ,y_n)=P_n(f;y_1, \ldots ,y_n)+Q,
\end{equation}
where $Q$ is either the zero matrix (if $y_j\neq z_1$), or $Q$ has
an $(n-1) \times (n-1)$ principal zero submatrix (if $y_j=z_1$ for some $j$),
and in any case $Q$ has at most one negative eigenvalue.
Let $k=\sqminus P_n(g;y_1, \ldots ,y_n)$.
Then (\ref{ag}) implies, using the Weyl inequalities for eigenvalues of the sum
of two Hermitian matrices (see, e.g., \cite[Section III.2]{B}) that
$\sqminus P_n(f;y_1, \ldots ,y_n)\leq k+1$.
\eproof
\medskip

{\bf Proof of $3\Rightarrow 1$ of Theorem \ref{T:1.4}.} Let $f$ be
defined on $\D\backslash \cQ$, where $\cQ$ is discrete. Furthermore,
let (\ref{2.10}) hold  for some integer $n\geq 0$ and let the set
$$
\cZ=\{z_1,\ldots,z_n\}\subset\D\backslash\cQ
$$
be such that
$$
\sqminus P_n(f; \, z_1,\ldots,z_n)=\kappa.
$$
Without loss of generality we can assume that
\begin{equation}
P_n(f; \, z_1,\ldots,z_n)=
\left[\frac{1-f(z_i)f(z_j)^*}{1-z_i{z}^*_j}\right]_{i,j=1}^n
\label{2.1}
\end{equation}
is not singular.
Indeed, if it happens that the matrix (\ref{2.1}) is singular, and its rank is
$m < n$, then by Lemma \ref{L:4.1}, there is a nonsingular  $m\times m$
principal submatrix $P_m(f;z_{i_1},\ldots,z_{i_m})$ of (\ref{2.1}).
A Schur complement argument shows that
$$
\sqminus P_m(f; \, z_{i_1},\ldots,z_{i_m})=\kappa.
$$
It follows then that ${\bf k}_m(f)=\kappa$. But then, in view of
(\ref{2.10}) and the nondecreasing property of the sequence
$\{{\bf k}_j(f)\}_{j=1}^{\infty}$,
we have ${\bf k}_{m+3}(f)=\kappa$, and the subsequent proof may be repeated
with $n$ replaced by $m$.

\smallskip

Setting $P_n(f; \, z_1,\ldots,z_n)=P$ for short, note the identity
\begin{equation}
P-TPT^*=FJF^*,
\label{2.2}
\end{equation}
where
\begin{equation}
T=\left[\begin{array}{ccc}{z}_1 && \\ &\ddots & \\
&& {z}_n\end{array}\right],\quad J=\left[\begin{array}{cr}1 &0\\
0&-1\end{array}\right],\quad
F=\left[\begin{array}{cc} 1 & f(z_1) \\
\vdots & \vdots \\ 1 & f(z_n)\end{array}\right].
\label{2.3}
\end{equation}
Consider the function
\begin{equation}
\Theta(z)=I_2-(1-z)F^*(I_n-zT^*)^{-1}P^{-1}(I_n-T)^{-1}FJ
\label{2.3a}
\end{equation}
It follows from (\ref{2.2}) by a straightforward standard
calculation (see, e.g., \cite[Section 7.1]{bgr}) that
\begin{equation}
J-\Theta(z)J\Theta(w)^*=(1-z{w}^*)F^*(I_n-zT^*)^{-1}P^{-1}
(I_n-{w}^*T)^{-1}F.
\label{2.3b}
\end{equation}
Note that $\Theta(z)$ is a rational function taking $J$--unitary values
on $\T$: $\Theta(z)J\Theta(z)^*=J$ for $z\in \T$. Therefore, by the symmetry
principle,
$$
\Theta(z)^{-1}=J\Theta(\frac{1}{{z}^*})^*J=
I_2+(1-z)F^*(I_n-T^*)^{-1}P^{-1}(zI_n-T)^{-1}FJ,
$$
which implies, in particular, that $\Theta(z)$ is invertible
at each point $z\not\in\cZ$.
By (\ref{2.10}) and by Lemma \ref{L:new},
\begin{equation}
\sqminus P_{n+3}(f; \, z_1,\ldots,z_n,\zeta_1,\zeta_2,\zeta_3)=\kappa
\label{2.4}
\end{equation}
for every choice of $\zeta_1, \, \zeta_2, \, \zeta_3\in\D\backslash\cQ$.
The matrix in (\ref{2.4}) can be written in the block form as
\begin{equation}
P_{n+3}(f; \, z_1,\ldots,z_n,\zeta_1,\zeta_2,\zeta_3)=
\left[\begin{array}{cc} P & \Psi^*\\ \Psi & P_3(f; \, \zeta_1, \zeta_2,
\zeta_3)\end{array}\right],
\label{2.6}
\end{equation}
where
$$
\Psi=\left[\begin{array}{c}\Psi_1 \\ \Psi_2 \\ \Psi_3 \end{array}\right]
\quad\mbox{and}\quad
\Psi_i=\left[\begin{array}{ccc}{\displaystyle
\frac{1-f(\zeta_i)f(z_1)^*}
{1-\zeta_i{z}^*_1}} &\ldots &
{\displaystyle\frac{1-f(\zeta_i)f(z_n)^*}{1-\zeta_i{z}^*_n}}
\end{array}\right] \; (i=1,2,3).
$$
The last formula for $\Psi_i$ can be written in terms of (\ref{2.3}) as
\begin{equation}
\Psi_i=[1 \; \; -f(\zeta_i)]F^*(I_n-\zeta_i T^*)^{-1}\qquad (i=1,2,3).
\label{2.5}
\end{equation}
By (\ref{2.10}), it follows from (\ref{2.6}) that
$$
P_3(f; \, \zeta_1, \zeta_2, \zeta_3)-\Psi P^{-1}\Psi^*\geq 0,
$$
or more explicitly,
\begin{equation}
\left[ \frac{1-f(\zeta_i)f(\zeta_j)^*}{1-\zeta_i{\zeta^*_j}}-
\Psi_i P^{-1}\Psi_j^*\right]_{i,j=1}^3\geq 0.
\label{2.8}
\end{equation}
By (\ref{2.5}) and (\ref{2.3b}),
\begin{eqnarray}
&&\frac{1-f(\zeta_i)f(\zeta_j)^*}{1-\zeta_i{\zeta^*_j}}-
\Psi_i P^{-1}\Psi_j^*\nonumber\\
&&=\left[\begin{array}{cc}1 &
-f(\zeta_i)\end{array}\right]\left\{\frac{J}{1-\zeta_i{\zeta}^*_j}-
F^*(I_n-\zeta_i
T^*)^{-1}P^{-1}(I_n-{\zeta}^*_jT)^{-1}F\right\}\left[\begin{array}{c}
1 \\ -f(\zeta_j)^*\end{array}\right]\nonumber \\
&&=\left[\begin{array}{cc}1 &
-f(\zeta_i)\end{array}\right]\frac{\Theta(\zeta_i)J\Theta(\zeta_j)^*}
{1-\zeta_i{\zeta}^*_j}\left[\begin{array}{c}
1 \\ -f(\zeta_j)^*\end{array}\right],\label{2.8a}
\end{eqnarray}
which allows us to rewrite (\ref{2.8}) as
\begin{equation}
\left[ \left[\begin{array}{cc}1 &
-f(\zeta_i)\end{array}\right]\frac{\Theta(\zeta_i)J\Theta(\zeta_j)^*}
{1-\zeta_i{\zeta}^*_j}\left[\begin{array}{c}
1 \\ -f(\zeta_j)^*\end{array}\right]\right]_{i,j=1}^3\geq 0.
\label{2.9}
\end{equation}
Introducing the block decomposition
$$
\Theta(z)=\left[\begin{array}{cc}\theta_{11}(z) & \theta_{12}(z)\\
\theta_{21}(z) & \theta_{22}(z)\end{array}\right]
$$
of $\Theta$ into four scalar blocks, note that
\begin{equation}
d(z):=\theta_{21}(z)f(z)-\theta_{11}(z)\neq 0, \quad
z\in \D\setminus (\cZ\cup\cQ).
\label{2.9a}
\end{equation}
Indeed, assuming that $\theta_{21}(\zeta)f(\zeta)=\theta_{11}(\zeta)$
for some $\zeta\in\D\backslash(\cZ\cup\cQ)$, we get
\begin{eqnarray*}
&&\left[\begin{array}{cc}1 &
-f(\zeta)\end{array}\right]\frac{\Theta(\zeta)J\Theta(\zeta)^*}
{1-|\zeta|^2}\left[\begin{array}{c}
1 \\ -f(\zeta)^*\end{array}\right]\nonumber\\
&&=-\left[\begin{array}{cc}1 &
-f(\zeta)\end{array}\right]\frac{\left[\begin{array}{c}\theta_{12}(\zeta)\\
\theta_{22}(\zeta)\end{array}\right]\left[\begin{array}{cc}
\theta_{12}(\zeta)^* & \theta_{22}(\zeta)^*\end{array}\right]}
{1-|\zeta|^2}\left[\begin{array}{c}1 \\ -f(\zeta)^*\end{array}\right]\leq
0,
\end{eqnarray*}
which contradicts (\ref{2.9}), unless $\det \Theta(\zeta)=0$. But as we
have mentioned above, $\Theta(z)$ is invertible at each point
$\zeta\not\in\cZ$.

\smallskip

Thus, the function
\begin{equation}
\sigma(z)=\frac{\theta_{12}(z)-f(z)\theta_{22}(z)}
{\theta_{21}(z)f(z)-\theta_{11}(z)}
\label{2.10z}
\end{equation}
is defined on $\D\backslash(\cZ\cup\cQ)$. Moreover,
$$
\left[\begin{array}{cc}1 & -f(\zeta)\end{array}\right]\Theta(\zeta)=-
d(\zeta)^{-1}\left[\begin{array}{cc}1
& -\sigma(\zeta)\end{array}\right]
$$
and thus, inequality (\ref{2.9}) can be written in terms of $\sigma$ as
$$
\left[ d(\zeta_i)^{-1}\left[\begin{array}{cc}1
& -\sigma(\zeta_i)\end{array}\right]\frac{J}
{1-\zeta_i{\zeta}^*_j}\left[\begin{array}{c}
1 \\ -\sigma(\zeta_j)^*\end{array}\right]
(d(\zeta_j)^*)^{-1}\right]_{i,j=1}^3\geq 0,
$$
or equivalently, (since $d(\zeta_i)\neq 0$) as
$$
\left[ \frac{1-\sigma(\zeta_i)\sigma(\zeta_j)^*}
{1-\zeta_i{\zeta}^*_j}\right]_{i,j=1}^3\geq 0.
$$
The latter inequality means that
$\; P_3(\sigma; \, \zeta_1,\zeta_2,\zeta_3)\geq 0 \;$
for every choice of points $\zeta_1, \, \zeta_2, \, \zeta_3$ in
$\D\backslash(\cZ\cup\cQ)$. By Corollary
\ref{T:1.2a}, $\sigma(z)$ is a Schur function. Although
$\sigma(z)$ has been defined via (\ref{2.10z}) on
$\D\backslash(\cZ\cup\{z_1,\ldots,z_n\})$, it admits an analytic
continuation to all of  $\D$, which still will be denoted by $\sigma(z)$.
It follows from (\ref{2.10z}) that $f$ coincides with the function
\begin{equation}
F(z)=\frac{\theta_{11}(z)\sigma(z)+\theta_{12}(z)}
{\theta_{21}(z)\sigma(z)+\theta_{22}(z)}.
\label{2.11}
\end{equation}
at every point $z\in\D\setminus(\cZ\cup\cQ)$.
Since $f$ has not been defined on $\cQ$, one can consider $F$
as a (unique) meromorphic extension of $f$. However,
$F$ need not coincide with $f$ on $\cZ$.

\smallskip

Now we prove that $f\in\cS_{\kappa}$. To this end it suffices to show that
\begin{equation}
\sqminus P_{n+r}(f; \, z_1,\ldots,z_n,\zeta_1,\ldots,\zeta_r)=\kappa
\label{2.13}
\end{equation}
for every choice of $\zeta_1,\ldots,\zeta_r\in\D\setminus (\cZ\cup\cQ)$.
Note that all possible ``jumps'' of $f$ are in $\cZ$ and at all other
points of $\D\setminus \cQ$, it holds that $f(\zeta)=F(\zeta)$. Thus,
writing
$$
P_{n+r}(f; \, z_1,\ldots,z_n,\zeta_1,\ldots,\zeta_r)=
\left[\begin{array}{cc} P & {\bf\Psi}^*\\ {\bf\Psi} & P_{r}(f; \,
\zeta_1, \ldots, \zeta_r)\end{array}\right],
$$
where
$$
{\bf\Psi}=\left[\begin{array}{c}\Psi_1 \\ \vdots \\ \Psi_r
\end{array}\right]
$$
and the $\Psi_i$'s are defined via (\ref{2.5}) for $i=1,\ldots,r$, we
conclude by the Schur complement argument that
\begin{equation}
\sqminus P_{n+r}(f; \, z_1,\ldots,z_n,\zeta_1,\ldots,\zeta_r)=
\sqminus P + \sqminus (P_{r}(f; \, ,\zeta_1, \ldots, \zeta_r)-\Psi
P^{-1}\Psi^*).
\label{2.14}
\end{equation}
It follows from the calculation (\ref{2.8a}) that
$$
P_{r}(f; \, \zeta_1, \ldots, \zeta_r)-\Psi
P^{-1}\Psi^*=\left[ \left[\begin{array}{cc}1 &
-f(\zeta_i)\end{array}\right]\frac{\Theta(\zeta_i)J\Theta(\zeta_j)^*}
{1-\zeta_i{\zeta}^*_j}\left[\begin{array}{c}
1 \\ -f(\zeta_j)^*\end{array}\right]\right]_{i,j=1}^r,
$$
which can be written in terms of functions $\sigma$ and $d$ defined in
(\ref{2.10z}) and (\ref{2.9a}), respectively, as
$$
P_{r}(f; \, \zeta_1, \ldots, \zeta_r)-\Psi
P^{-1}\Psi^*=\left[
d(\zeta_i)^{-1}\frac{1-\sigma(\zeta_i)\sigma(\zeta_j)^*}
{1-\zeta_i{\zeta}^*_j}(d(\zeta_j)^*)^{-1}\right]_{i,j=1}^r.
$$
We have already proved that $\sigma$ is a Schur function and therefore,
$$
P_{r}(f; \, \zeta_1, \ldots, \zeta_r)-\Psi
P^{-1}\Psi^*\geq 0,
$$
which together with (\ref{2.14}) implies (\ref{2.13}).
\qed

\section{Proof of Theorem \ref{T:3.1}}
\setcounter{equation}{0}

Since the first inequality in
(\ref{3.1}) is obvious, we need to prove only the second
inequality.
\smallskip

If $\tilde{f}$ is the meromorphic part of $f$, then
by Theorem \ref{T:1.2b} $\tilde{f}\in\cS_{q}$, and therefore
by Lemma \ref{L:last1}, ${\bf k}_n(f)\leq q+\ell,$ for $n=1,2, \ldots $.
Therefore, it suffices to prove that there exist $q+2\ell$ distinct
points
$u_1, \ldots , u_{q+2\ell}\in{\rm Dom}\,(f)$ such that
$$
\sqminus P_{q + 2\ell}(f;u_1, \ldots ,u_{q  +2\ell})\geq q+\ell.
$$

We start with notation and some preliminary
results. Let
$$
J_r(a)=\left[\begin{array}{ccccc} a &0 & 0 & \cdots &0 \\
1& a & 0 & \cdots &0 \\ 0&1&a& \cdots &0 \\
\vdots & \vdots & \vdots & \ddots & \vdots \\ 0&0&0& \cdots 1 & a
\end{array}\right], \qquad a\in \C
$$
be the lower triangular $r \times r$ Jordan block with
eigenvalue $a$ and let $E_r$ and $G_r$ be vectors from $\C^r$ defined by
\begin{equation}
\label{6.01}
E_r=\left[\begin{array}{c} 1 \\ 0 \\ \vdots \\ 0\end{array}\right]\in\C^r
\quad \mbox{and}\quad
G_r=\left[\begin{array}{c} 1 \\ 1 \\ \vdots \\ 1\end{array}\right]\in\C^r.
\end{equation}
Given an ordered set $\cZ=\{z_1, \ldots ,z_k\}$ of
distinct points in the complex plane, we denote by $\Phi(\cZ)$ the lower
triangular $k \times k$ matrix defined by
\begin{equation}\label{x.1}
\Phi(\cZ)=\left[\Phi_{i,j}\right]_{i,j=1}^k, \qquad
\Phi_{i,j}=\left\{ \begin{array}{cl}
{\displaystyle\frac{1}{\phi_i^{\prime} (z_j)}} &
{\rm if} \ i\geq j, \\ 0 & {\rm if} \ i<j, \end{array}\right.
\end{equation}
where $\phi_i(z)=\prod_{j=1}^i (z-z_j)$. Furthermore, given a set $\cZ$ as
above, for a complex valued function $v(z)$ we define recursively the divided
differences $[z_1, \ldots ,z_j]_v$ by
$$
[z_1]_v=v(z_1), \qquad [z_1, \ldots ,z_{j+1}]_v=
\frac{[z_1, \ldots ,z_{j}]_v-[z_2, \ldots ,z_{j+1}]_v}{z_1- z_{j+1}}, \qquad
j=1, \ldots ,k,
$$
and use the following notation for associated matrices and vectors:
$$ D(\cZ)=
\left[\begin{array}{cccc} z_1 &0 & \cdots &0 \\
0& z_2  & \cdots &0 \\
\vdots & \vdots  & \ddots & \vdots \\ 0&0& \cdots  & z_k
\end{array}\right], \quad J(\cZ)=
\left[\begin{array}{ccccc} z_1 &0 & 0 & \cdots &0 \\
1& z_2 & 0 & \cdots &0 \\ 0&1&z_3& \cdots &0 \\
\vdots & \vdots & \vdots & \ddots & \vdots \\ 0&0&0& \cdots 1 & z_k
\end{array}\right],
$$
\begin{equation}
v(\cZ)=\left[\begin{array}{c} v(z_1) \\
v(z_2) \\ \vdots \\
v(z_k)\end{array}\right], \quad
[\cZ]_v= \left[\begin{array}{c} [z_1]_v \\ \left[z_1, z_2\right]_v \\
\vdots \\ \left[z_1, \ldots ,z_k\right]_v \end{array}\right].
\label{x.2.2a}
\end{equation}
\begin{La}
\label{L:x.1}
Let  $\cZ=\{z_1, \ldots,z_k\}$ be an ordered set of distinct points,
and let $\Phi(\cZ)$ be defined as in $(\ref{x.1})$. Then:
\begin{equation}
\Phi(\cZ)D(\cZ)=J(\cZ)\Phi(\cZ) \qquad {\rm and } \qquad
\Phi(\cZ)v(\cZ)=[\cZ]_v.
\label{x.2.3}
\end{equation}
Moreover, if the function $v(z)$ is analytic at
$z_0\in\C$, with the Taylor series $v(z)=\sum_{k=0}^{\infty} v_k(z-z_0)^k$,
then
\begin{equation}
\label{x.4}
\lim_{z_1,\ldots,z_k\rightarrow z_0}
\Phi(\cZ)v(\cZ)=\left[\begin{array}{c} v_0\\ v_1 \\ \vdots \\ v_{k-1}
\end{array}\right]. \end{equation}
\end{La}
{\bf Proof:} Formulas (\ref{x.2.3}) are verified by direct computation;
formula (\ref{x.4}) follows from (\ref{x.2.3}), since
$$
\lim_{z_i\rightarrow z_0,\ i=1,\ldots ,j} [z_1, \ldots ,z_j]_v=
\frac{v^{(j-1)}(z_0)}{(j-1)!}=v_{j-1}, \quad  j=1, \ldots ,k .
$$
In the sequel it will be convenient to use the following notation:
${\rm diag}\, \left(X_1, \ldots ,X_k\right)$ stands for the block diagonal
matrix with the diagonal blocks $X_1, \ldots ,X_k$ (in that order).
\begin{La}
\label{L:x.2}
Let $w_1, \ldots ,w_k$ be distinct points in $\D$, and let $K$ be the unique
solution of the Stein equation
\begin{equation}\label{x.4a}
K -AKA^*=EE^*, \end{equation}
where
\begin{equation}
A=\left[\begin{array}{ccc}J_{r_1}(w_1) && \\ & \ddots &  \\ &&
J_{r_k}(w_k)\end{array}\right], \quad E=\left[\begin{array}{c}E_{r_1} \\
\vdots \\ E_{r_k}\end{array}\right]
\label{x.40}
\end{equation}
and $E_{r_j}$ are defined via the first formula in $(\ref{6.01})$.
Then the normalized Blaschke product
\begin{equation}\label{x.4b}
b(z)=e^{i\alpha} \prod_{j=1}^k \left(\frac{z-w_j}{1-zw_j^*}\right)^{r_j},\qquad
b(1)=1,
\end{equation}
admits a realization
\begin{equation}\label{x.4c}
b(z)=1+(z-1)E^*(I-zA^*)^{-1}K^{-1}(I-A)^{-1}E,\end{equation}
and the following formula holds:
\begin{equation}\label{x.4d}
1-b(z)b(w)^*=(1-zw^*)E^*(I-zA^*)^{-1}K^{-1}(I-w^*A)^{-1}E, \quad z,w\in \D.
\end{equation}
\end{La}
{\bf Proof:} First we note that
$K$ is given by the convergent series
\begin{equation}\label{x.4e}
K=\sum_{j=0}^{\infty} A^jEE^*(A^*)^j.
\end{equation}
Since the pair $(E^*, A^*)$ is observable, i.e, $\cap_{j=0}^{\infty}
{\rm Ker}\, (E^*(A^*)^j)=\{0\}$, the matrix $K$ is positive definite. Equality
(\ref{x.4d}) follows from (\ref{x.4c}) and (\ref{x.4a}) by a standard
straightforward  calculation (or from (\ref{2.3b}) upon setting $J=I_2$
in (\ref{2.3})--(\ref{2.3b})).
It follows from (\ref{x.4d}) that the function $b(z)$ defined via
(\ref{x.4c}) is inner.
Using the fact that ${\rm det}\, (I+XY)={\rm det}\,(I+YX)$ for
matrices $X$ and $Y$ of sizes $u \times v$ and $v \times u$ respectively,
we obtain from (\ref{x.4a}) and (\ref{x.4c}):
\begin{eqnarray*}
b(z)&= & {\rm det}\, \left(
I+(z-1)(I-zA^*)^{-1}K^{-1}(I-A)^{-1}(K-AKA^*)\right) \\
&=& {\rm det}\, \left((I-zA^*)^{-1}K^{-1}(I-A)^{-1}\right)
 \\ && \times \  {\rm det} \, \left((I-A)K(I-zA^*)+(z-1)(K-AKA^*)\right) \\
&= & {\rm det}\, \left((I-zA^*)^{-1}K^{-1}(I-A)^{-1}\right)\, \cdot\,
{\rm det}\, \left((zI-A)K(I-A^*)\right) \\
&= & c\, \frac{{\rm det}\,(zI -A)}{{\rm det}\, (I-zA^*)} \qquad
{\rm for}\ \ {\rm some} \ \ c\in \C, \ |c|=1.
\end{eqnarray*}
It follows that the degree of $b(z)$ is equal to $r:=\sum_{j=1}^k
r_j$, and $b(z)$ has zeros at $w_1, \ldots ,w_k$ of multiplicities
$r_1, \ldots ,r_k$, respectively.
Since $b(1)=1$,
the function $b(z)$ indeed
coincides with (\ref{x.4b}). \qed
\bigskip

The result of Lemma \ref{L:x.2} is known also for matrix valued inner functions
(see
\cite[Section 7.4]{bgr}, where it is given in a slightly different form for
$J$-unitary matrix functions), in which case it can be interpreted as a formula
for
such functions having a prescribed left null pair with respect to the unit disk
$\D$ (see \cite{bgr} and relevant references there).

Let
\begin{equation}
\label{x.5}
f(z)=\left\{\begin{array}{cl} {\displaystyle \frac{S(z)}{b(z)}} &
{\rm if}\, z\not\in \{z_1,\ldots,z_{\ell}, w_{t+1},\ldots ,w_k\} \\
f_j & {\rm if }\, z=z_j, j=1, \ldots , \ell \end{array}\right.
\end{equation}
where $S(z)$ is a Schur function not vanishing at any of the $z_j$'s, $b(z)$
is the Blaschke product given by (\ref{x.4b}). The points $z_1, \ldots
,z_{\ell}$ are assumed to be distinct points in $\D$, and $w_1,\ldots ,w_k$ are
also assumed to be distinct in $\D$. Furthermore, we assume that
$w_j=z_j$ for  $j=1, \ldots , t$, $\{w_{t+1}, \ldots ,w_k\}\cap
\{z_{t+1},\ldots
,z_{\ell}\}=\emptyset$, and ${\displaystyle \frac{S(z_j)}{b(z_j)}}\neq f_j$ for
$j=t+1, \ldots ,\ell.
$ The cases when $t=0$, i.e.,
$\{w_{1}, \ldots ,w_k\}\cap \{z_{1},\ldots
,z_{\ell}\}=\emptyset$, and when $t=\min\{k,\ell\}$ are not excluded; in these
cases the subsequent arguments should be modified in obvious ways.
We let $N=2\ell+\sum_{j=1}^kr_j$.

Take $N$ distinct points in the unit disk sorted into $k+2$ ordered sets
\begin{equation}\label{x.5a}
 \cM_j=\{\mu_{j,1},\ldots ,\mu_{j,r_j}\}, \ j=1, \ldots, k; \quad
\cN=\{\nu_1, \ldots ,\nu_{\ell} \}; \quad \cZ=\{z_1,\ldots ,z_{\ell}\}
\end{equation}
and such that $\cM_j\cap \cW=\emptyset$, $j=1, \ldots ,k$, and
$\cN\cap \cW=\emptyset$, where $\cW=\{w_1, \ldots , w_k\}$. Consider the
corresponding Pick matrix
$$ P=P_N(f;\mu_{1,1},\ldots , \mu_{k,r_k}, \nu_1, \ldots , \nu_{\ell},z_1,
\ldots , z_{\ell}). $$
We shall show that if $\mu_{j,i}$ and $\nu_j$ are sufficienly close to
$w_j$ and $z_{j}$, respectively, then $\sqminus P\geq N-\ell$.

It is readily seen that $P$ is a unique solution of the Stein equation
\begin{equation}
\label{x.6}
P-TPT^*=G_NG_N^*-CC^* \end{equation}
where $G_N\in\C^N$ is defined via the second formula in (\ref{6.01}) and
$$
T=\left[\begin{array}{ccccc}D(\cM_1) &&&&\\ &\ddots &&& \\
&&D(\cM_k)&&\\ &&& D(\cN) &\\ &&&& D(\cZ)\end{array}\right],\qquad
C= \left[\begin{array}{c} f(\cM_1) \\ \vdots  \\ f(\cM_k) \\ f(\cN) \\ f(\cZ)
\end{array}\right].
$$
We recall that by definition (\ref{x.2.2a}) and in view of (\ref{x.5}),
$$
f(\cM_j)=\left[\begin{array}{c}f(\mu_{j,1}) \\ \vdots \\ f(\mu_{j,r_j})
\end{array}\right],\quad
f(\cN)=\left[\begin{array}{c}f(\nu_1)\\ \vdots \\
f(\nu_\ell)\end{array}\right],\quad
f(\cZ)=\left[\begin{array}{c}f_1\\ \vdots \\ f_{\ell}\end{array}\right].
$$
Consider the matrices
$$
B_j=\Phi(\cM_j){\rm diag}\,\left(b(\mu_{j,1}), \ldots ,
b(\mu_{j,r_j})\right)\quad (j=1,\ldots,k)
$$
and note that by Lemma \ref{L:x.1},
\begin{eqnarray}
B_jD(\cM_j)&=&\Phi(\cM_j)D(\cM_j){\rm diag}\,\left(b(\mu_{j,1}), \ldots ,
b(\mu_{j,r_j})\right)\nonumber\\
&=&J(\cM_j)\Phi(\cM_j){\rm diag}\,\left(b(\mu_{j,1}),\ldots,b(\mu_{j,r_j})
\right)\nonumber\\ &=&J(\cM_j)B_j,\nonumber\\
B_jf(\cM_j)&=&\Phi(\cM_j)S(\cM_j)=\left[\cM_j\right]_S,\nonumber\\
B_jG_{r_j}&=&\Phi(\cM_j)b(\cM_j)=\left[\cM_j\right]_b.\nonumber
\end{eqnarray}
The three last equalities together with block structure of the
matrices $T$, $C$ and $G_N$, lead to
\begin{equation}
\label{x.7}
BT=T_1B, \quad Y_1=BY, \quad {\rm and} \quad C_1=BC,
\end{equation}
where
\begin{equation}
\label{x.6ab}
B={\rm diag}\,\left(B_1,\ldots,B_j, \, I_{2\ell}\right), \quad
T_1={\rm diag}\, \left(J(\cM_1), \ldots,
J(\cM_k), D(\cN), D(\cZ)\right), \end{equation}
\begin{equation} \label{x.6ac} Y_1=\left[\begin{array}{c}
\left[\cM_1\right]_b \\ \vdots \\
\left[\cM_k\right]_b \\ G_{2\ell}\end{array}\right]\quad\mbox{and}\quad
C_1 \ = \  \left[\begin{array}{c} \left[\cM_1\right]_S \\
 \vdots \\
\left[\cM_k\right]_S \\ f(\cN) \\ f(\cZ) \end{array}\right].
\end{equation}
Pre- and post-multiplying (\ref{x.6}) by
$B$ and $B^*$, respectively, we conclude, on account of (\ref{x.7}), that the
matrix $
P_1:=BPB^* $
is the unique solution of the Stein equation
\begin{equation}\label{x.9}
P_1-T_1P_1T_1^*=Y_1Y_1^*-C_1C_1^*,
\end{equation}
where $T_1$, $Y_1$, and $C_1$ are given by (\ref{x.6ab}) and
(\ref{x.6ac}).

\smallskip

Recall that all the entries in (\ref{x.9}) depend on the $\mu_{j,i}$'s.
We now let $\mu_{j,i} \rightarrow w_j$, for $i=1, \ldots ,r_j$,
$j=1,\ldots , k$. Since $b$ has zero of order $r_j$ at $w_j$, it follows
 by (\ref{x.4}) that
$$
\lim_{\mu_{j,i}\rightarrow w_j}[\cM_j]_b=0, \quad j=1, \ldots, k,
$$
and thus,
$$
Y_2:=\lim_{\mu_{j,i}\rightarrow w_j} Y_1=\left[\begin{array}{c}
0 \\ G_{2\ell}\end{array}\right]\in\C^{N}.
$$
Similarly, we get by (\ref{x.2.2a}),
$$
C_2:=\lim_{\mu_{j,i}\rightarrow w_j} C_1=\left[\begin{array}{c}
\widehat{S}_1 \\ \vdots \\ \widehat{S}_k \\ f(\cN) \\ f(\cZ)
\end{array}\right], \quad\mbox{where}\quad
\widehat{S}_j=\left[\begin{array}{c} S(w_j) \\[1mm]
{\displaystyle \frac{S'(w_j)}{1!}}
\\[1mm]
\vdots \\[1mm]
{\displaystyle \frac{S^{(r_j-1)}(w_j)}{(r_j-1)!}}\end{array}\right]
\equiv \left[\begin{array}{c}s_{j,0} \\ s_{j,1} \\
\vdots \\ s_{j,r_{j-1}}\end{array}\right].
$$
Furthermore, taking the limit in (\ref{x.6ab}) as $\mu_{j,i} \rightarrow
w_j$, we get
\begin{eqnarray*}
T_2:=\lim_{\mu_{j,i}\rightarrow w_j} T_1 &=&{\rm diag} \,
\left(J_{r_1}(w_1),
\ldots , J_{r_k}(w_k),D(\cN),D(\cZ) \right) \nonumber\\
&=& {\rm diag}\, (A, D(\cN), D(\cZ)),
\end{eqnarray*}
where $A$ is given in (\ref{x.40}).
Since the above three limits exist, there also exists the limit
\begin{equation}
\label{x.13}
P_2:= \lim_{\mu_{j,i}\rightarrow w_j} P_1,
\end{equation}
which serves to define a unique solution $P_2$ of the Stein equation
\begin{equation}\label{x.14}
P_2-T_2P_2T_2^*=Y_2Y_2^* - C_2C_2^*.
\end{equation}
Let $\cS_j$ be the lower triangular Toeplitz matrices defined by:
$$
\cS_j=\left[\begin{array}{ccccc} s_{j,0} & 0& \ldots &0&0 \\
s_{j,1} & s_{j,0} & \ldots &0&0 \\ \vdots & \ddots & \ddots & \vdots & \vdots
\\ s_{j,r_j-2} &  &\ldots & s_{j,0}&0 \\ s_{j,r_j-1} & s_{j,r_j-2} & \ldots
& s_{j,1} & s_{j,0} \end{array}\right], \qquad j=1, \ldots ,k,
$$
which are invertible because $S(w_j)=s_{j,0}\neq 0$, $j=1,\ldots, k$. Let
$$
R:={\rm diag}\,(\cS_1^{-1}, \ldots , \, \cS_{k}^{-1}, \,
I_{2\ell})\quad\mbox{and}\quad
C_3=\left[\begin{array}{c}E \\ f(\cN)\\f(\cZ)\end{array}\right],
$$
where $E$ is given in (\ref{x.40}). The block structure of matrices
$R$, $T_2$, $C_2$, $C_3$, $E$ and $Y_2$ together with selfevident
relations
$$
\cS_j^{-1}J_{r_j}(w_j)=J_{r_j}(w_j)\cS_j^{-1},\quad
\cS_j^{-1}\widehat{S}_1=E_{r_j}\quad (j=1,\ldots,k),
$$
lead to
$$
RT_2=T_2R,\quad \quad RC_2=C_3\quad {\rm and }\quad RY_2=Y_2.
$$
Taking into account the three last equalities, we pre- and post-multiply
(\ref{x.14}) by $R$ and $R^*$, respectively, to conclude that the matrix
$P_3:=RP_2R^*$ is the unique solution of the Stein equation
\begin{equation}
\label{x.16}
P_3-T_2P_3T_2^*=Y_2Y_2^*-C_3C_3^*.
\end{equation}
Denote
$$
\cN_1=\{\nu_1,\ldots ,\nu_t\}, \qquad \cZ_1 =\{z_1, \ldots ,z_t\},
$$ $$
 \cN_2=\{\nu_{t+1},\ldots ,\nu_\ell\}, \qquad  \cZ_2 =\{z_{t+1}, \ldots
,z_\ell\}.
$$
By the hypothesis, $\cZ_2 \cap \cW=\emptyset.$ Let $\nu_j \rightarrow
z_j$ for $j=t+1, \ldots , \ell$. Then $$f(\cN_2)\rightarrow \frac{S}{b}(\cZ_2),
$$ and note that by the assumptions on $f$,
\begin{equation}\label{x.16a}
\frac{S}{b} (z_j)\neq f(z_j), \qquad j=t+1, \ldots ,\ell. \end{equation}
Taking limits in (\ref{x.16}) as $\nu_j \rightarrow z_j$ ($j=t+1,\ldots
,\ell$), we arrive at the equality
\begin{equation}\label{x.17}
P_4- T_3P_4T_3^*=Y_2Y_2^*- C_4C_4^*,
\qquad P_4:=\lim_{\nu_j\rightarrow z_j,\ j=t+1,\ldots ,\ell} P_3,
\end{equation}
where
$$
T_3={\rm diag}\, \left( A, D(\cN_1), D(\cZ_2), D(\cZ_1), D(\cZ_2) \right),
\quad C_4 = \left[\begin{array}{c} E \\  f(\cN_1)\\ {\displaystyle\frac{S}{b}
(\cZ_2)}\\
f(\cZ_1) \\
f(\cZ_2)\end{array} \right].
$$
By (\ref{x.16a}), the matrix
$$
L:={\rm diag} \, \left( \left(\frac{S}{b} - f\right)(z_{t+1}), \ldots ,
\left(\frac{S}{b}-f\right)(z_{\ell}) \right)
$$
is invertible. Let
\begin{equation}\label{x.18a}
T_4=\left[\begin{array}{cccc} A &0 &0 &0 \\ 0& D(\cZ_2)&0&0\\
0&0& D(\cZ_1)&0 \\ 0&0&0& D(\cN_1)\end{array}\right],\quad
C_5= \left[\begin{array}{c} E \\ G_{\ell-t} \\ f(\cZ_1) \\ f(\cN_1)
\end{array}\right],\quad
Y_3= \left[\begin{array}{c} 0 \\ G_{2t}\end{array}\right]
\end{equation}
and
$$
X=\left[\begin{array}{ccccc} I_{N-2\ell} & 0 & 0 & 0 & 0\\
0 & 0 & L^{-1} & 0& -L^{-1} \\
0 & 0 & 0& I_t & 0 \\ 0 & I_t &0 & 0&0 \end{array}\right].
$$
Since
$$
XT_3=T_4X, \quad XC_4=C_5,\quad XY_2=Y_3,
$$
pre- and post-multiplication of (\ref{x.17}) by
$X$ and by $X^*$, respectively,
leads to the conclusion that the matrix $\; P_5:=XP_4X^* \;$
is the unique solution of the Stein equation
\begin{equation}\label{x.19}
P_5-T_4P_5T_4^*=Y_3Y_3^*-C_5C_5^*. \end{equation}
By (\ref{x.18a}) and (\ref{x.19}), $P_5$ has
the form
$$ P_5=\left[\begin{array}{c|ccc} -K & Q_1^* & Q_2^*& Q_3^* \\
\hline Q_1 & R_{11}&R_{12} & R_{13}\\ Q_2 & R_{21}& R_{22} & R_{23}\\
Q_3 & R_{31}& R_{32}& R_{33}
\end{array}\right], $$ where
$K$ is the matrix given by (\ref{x.4e}), and $Q_j$, $j=1,2,3$, are solutions
of
\begin{eqnarray*}
Q_1-D(\cZ_2)Q_1A^*&= &-G_{\ell-t}E^*, \\
Q_2-D(\cZ_1)Q_2A^*&=&-f(\cZ_1)E^*,
\\
 Q_3-D(\cN_1)Q_3A^*&=& -f(\cN_1)E^*, \end{eqnarray*}
respectively. Thus, denoting by $\left[Q_{\alpha}\right]_j$ the $j$th row of
$Q_{\alpha}$, $\alpha=1,2,3$, we have from the three last equalities (recall
that
$w_j=z_j$ for $j=1, \ldots ,t$):
\begin{eqnarray*}
\left[Q_1\right]_j&=& -E^*(I-z_{j+t}A^*)^{-1}, \qquad j=1, \ldots, \ell-t; \\
\left[Q_2\right]_j&=& -f_jE^*(I-w_{j}A^*)^{-1}, \qquad j=1, \ldots, t; \\
\left[Q_3\right]_j&=& -f(\nu_j)E^*(I-\nu_jA^*)^{-1}, \qquad j=1, \ldots,
t. \end{eqnarray*}
Furthermore, in view of the three last relations and (\ref{x.4d}), and
since $b(w_j)=0$, following equalities are obtained:
\begin{eqnarray*}
\left[Q_1\right]_jK^{-1}\left[Q_1\right]_i^*&= &
\frac{1-b(z_{t+j})b(z_{t+i})^*}{1-z_{t+j}z_{t+i}^*} , \quad j,i=1, \ldots
,\ell-t,\\
\left[Q_2\right]_jK^{-1}\left[Q_1\right]_i^*&= &
\frac{f_j}{1-w_{j}z_{t+i}^*} , \quad j=1, \ldots ,t; \ i=1, \ldots ,\ell-t,\\
\left[Q_3\right]_jK^{-1}\left[Q_1\right]_i^*&= &
f(\nu_j)\,\frac{1-b(\nu_{j})b(z_{t+i})^*}{1-\nu_{j}z_{t+i}^*} ,
\quad j=1, \ldots ,t; \ i=1, \ldots ,\ell-t,\\
\left[Q_2\right]_jK^{-1}\left[Q_2\right]_i^*&= &
\frac{f_jf_i^*}{1-w_{j}w_{i}^*} , \quad j,i=1, \ldots ,t, \\
\left[Q_3\right]_jK^{-1}\left[Q_2\right]_i^*&= &
\frac{f(\nu_j)f_i^*}{1-\nu_{j}w_{i}^*} , \quad j,i=1, \ldots ,t,\\
\left[Q_3\right]_jK^{-1}\left[Q_3\right]_i^*&= &
f(\nu_j)f(\nu_i)^*\frac{1-b(\nu_{j})b(\nu_{i})^*}{1-\nu_{j}\nu_{i}^*},
\quad j,i=1, \ldots ,t.
\end{eqnarray*}
Next, the Schur complements are computed:
\begin{eqnarray*}
\left[R_{11}+Q_1K^{-1}Q_1^*\right]_{ji}&=&
\frac{-1}{1-z_{t+j}z_{t+i}^*}+
\frac{1-b(z_{t+j})b(z_{t+i})^*}{1-z_{t+j}z_{t+i}^*}\\
&=&-\frac{b(z_{t+j})
b(z_{t+i})^*}{1-z_{t+j}z_{t+i}^*}, \\
\left[R_{21}+Q_2K^{-1}Q_1^*\right]_{ji}&= &0 , \\
\left[ R_{31}+Q_3K^{-1}Q_1^*\right]_{ji} & = &
-\frac{f(\nu_j)}{1-\nu_jz_{t+i}^*}+f(\nu_j)
\frac{1-b(\nu_j)b(z_{t+i})^*}{1-\nu_jz_{t+i}^*} \\
&= & -f(\nu_j)\frac{b(\nu_j)b(z_{t+i})^*}{1-\nu_jz_{t+i}^*}\\
&=& -\frac{S(\nu_j)b(z_{t+i})^*}{1-\nu_jz_{t+i}^*},
\end{eqnarray*} \begin{eqnarray*}
\left[R_{22}+Q_2K^{-1}Q_2^*\right]_{ji}
&=& \frac{1-f_jf_i^*}{1-w_jw_i^*}+
 \frac{f_jf_i^*}{1-w_jw_i^*} =  \frac{1}{1-w_jw_i^*}, \\
\left[R_{32}+Q_3K^{-1}Q_2^*\right]_{ji}&=&
\frac{1-f(\nu_j)f_i^*}{1-\nu_jw_i^*}+
 \frac{f(\nu_j)f_i^*}{1-\nu_jw_i^*} =  \frac{1}{1-\nu_jw_i^*}, \\
\left[ R_{33}+Q_3K^{-1}Q_3^*\right]_{ji} & = &
\frac{1-f(\nu_j)f(\nu_i)^*}{1-\nu_j\nu_i^*}+f(\nu_j)f(\nu_i)^*
\frac{1-b(\nu_j)b(\nu_i)^*}{1-\nu_j\nu_i^*} \\
&= & \frac{1-S(\nu_j)S(\nu_i)^*}{1-\nu_j\nu_i^*}\, .
\end{eqnarray*}
In the obtained Schur complement $\left[R_{i,j}+Q_iK^{-1}Q_j
\right]_{i,j=1}^3$, we let $\nu_j\rightarrow w_j$ ($j=1,2,\ldots ,t$) and
then pre- and  post-multiply by the matrix $\left[\begin{array}{cc}
b(\cZ_2)^{-1} & 0 \\  0 & I_{2t}\end{array}\right]$ and by its adjoint,
respectively, resulting in the matrix
\begin{equation}
\label{x.20}
\left[\begin{array}{ccc} -K_{11}&0 & K_{31}^*\\
0 & K_{22} & K_{22}\\ K_{31} & K_{22} &
\left[{\displaystyle \frac{1-S(w_j)S(w_i)^*}{1-w_jw_i^*}}\right]_{j,i=1}^t
\end{array}\right],
\end{equation}
where
\begin{equation}\label{x.20a}
K_{11}= \left[\frac{1}{1-z_{t+j}z_{t+i}^*}
\right]_{j,i=1}^{\ell-t}, \ \
K_{22}=\left[\frac{1}{1-w_jw_i^*}\right]_{j,i=1}^t,
\ \ K_{31}=\left[-\frac{S(w_j)}{1-w_jz_{t+i}^*}\right]_{j,i=1}^{t,\ell-t}.
\end{equation}
Note that the $j$-th row $[K_{31}]_j$ of $K_{31}$ can be written in the
form
$$
[K_{31}]_j=S(w_j)G_{\ell-t}^*(I-w_jD(\cZ_2)^*)^{-1}\quad (j=1,\ldots,t).
$$
Using Lemma \ref{L:new}, in view of
the reductions made so far from $P$ to $P_5$, to prove that $\sqminus
P\geq N-\ell$, we only have to show that
the Schur complement ${\bf S}$ to the block
$\left[\begin{array}{cc} -K_{11}&0\\ 0 &K_{22} \end{array}\right]$ in
(\ref{x.20}) has at least $t$ negative eigenvalues, i.e., is negative
definite. This Schur complement is equal to
\begin{eqnarray}
{\bf S}&=&\left[\frac{1-S(w_j)S(w_i)^*}{1-w_jw_i^*}\right]_{j,i=1}^t
+K_{31}K_{11}^{-1}K_{31}^* -K_{22}\nonumber\\
&=&-\left[\frac{S(w_j)S(w_i)^*}{1-w_jw_i^*}\right]_{j,i=1}^t
+K_{31}K_{11}^{-1}K_{31}^*\nonumber\\
&=&-\left[\frac{S(w_j)S(w_i)^*}{1-w_jw_i^*}-[K_{31}]_jK_{11}^{-1}
\left([K_{31}]_i\right)^*\right]_{j,i=1}^t\nonumber \\
&=&-\left[\frac{S(w_j)S(w_i)^*}{1-w_jw_i^*}\right.\nonumber\\
&&\quad \left. -S(w_j)G_{\ell-t}(I-w_jD(\cZ_2)^*)^{-1}K_{11}^{-1}
(I-w_i^*D(\cZ_2))^{-1}G_{\ell-t}^*S(w_i)^*\right]_{j,i=1}^t.
\label{x.22}
\end{eqnarray}
Introduce the rational function
\begin{equation}
\label{x.20b}
\vartheta(z)=1+(z-1)G_{\ell-t}^*(I-zD(\cZ_2)^*)^{-1}K_{11}^{-1}
(I-D(\cZ_2))^{-1}G_{\ell-t}.
\end{equation}
and note that $K_{11}$ satisfies the Stein equation
$$
K_{11}-D(\cZ_2)K_{11}D(\cZ_2)^*=G_{\ell-t}G_{\ell-t}^*.
$$
Upon setting $k=\ell-t$ and $r_1=\ldots=r_k=1$ and
picking points $z_{t+1},\ldots,z_\ell$ instead of $w_1,\ldots,w_k$
in Lemma \ref{L:x.2}, we get $A=D(\cZ_2)$, $K=K_{11}$, $E=G_{\ell-t}$
and thus, by Lemma
\ref{L:x.2}, we conclude that
\begin{equation}
\label{x.22c}
\vartheta(z)=e^{i\alpha} \prod_{j=1}^{\ell-t}
\frac{z-z_{t+j}}{1-zz_{t+j}^*},
\end{equation}
where $\alpha\in\R$ is chosen to provide the normalization
$\vartheta(1)=1$, and moreover, relation (\ref{x.4d}) in the present
setting takes the form
\begin{equation}
\label{x.22a}
1-\vartheta(z)\vartheta(w)^*=(1-zw^*)G_{\ell-t}^*(I-zD(\cZ_2)^*)^{-1}
K_{11}^{-1}(I-w^*D(\cZ_2))^{-1}G_{\ell-t}.
\end{equation}
It follows from (\ref{x.22c}) that $Z(\vartheta)=\cZ_2$ (more precisely,
$\vartheta$ is of degree $\ell-t$ and has simple poles at $z_{t+1},
\ldots,z_\ell$). Since $\cZ_1\cap\cZ_2=\emptyset$, it follows that
$\vartheta(w_j)\neq 0$, $j=1, \ldots, t$. Upon making use of
(\ref{x.22a}), we rewrite (\ref{x.22}) as
\begin{equation}
\label{x.23}
{\bf S}= -\left[\frac{S(w_j)\vartheta(w_j)\vartheta(w_i)^*S(w_i)^*}{1
-w_jw_i^*}\right]_{j,i=1}^t
\end{equation}
and conclude, since $S(w_j)\vartheta(w_j)\neq 0$ ($j=1, \ldots, t$)
that ${\bf S}$ is negative definite. Summarizing, we have
\begin{eqnarray*}
\sqminus P\geq \sqminus P_5&=&\sqminus (-K)+\sqminus
\left[\begin{array}{cc} -K_{11}&0\\ 0 &K_{22} \end{array}\right]
+\sqminus {\bf S}\\ &=&(N-2\ell)+(\ell-t)+t=N-\ell,
\end{eqnarray*}
which completes the proof.\qed

\section{Theorems \ref{T:1.4} and \ref{T:1.new}: proofs}
\setcounter{equation}{0}

{\bf Proof of Theorem \ref{T:1.4}:}
The implication {\bf 3} $\Rightarrow$ {\bf 1} was already proved,
and
{\bf 1} $\Rightarrow$ {\bf 3} follows from the definition of $\cS_{\kappa}$.
\smallskip

Assume {\bf 2} holds, and let $\tilde{f}$ be the standard function that
extends $f$ as stated in {\bf 2}. By Theorem \ref{T:3.1} $\tilde{f}\in
\cS_{\kappa}$. It is obvious from the definition of $\cS_{\kappa}$ that
we have $f\in \cS_{\kappa'}$ for some $\kappa'\leq\kappa$. However,
Remark \ref{R:last} implies that in fact $\kappa'=\kappa$, and {\bf 1}
follows.
\smallskip

It remains to prove {\bf 3} $\Rightarrow $ {\bf 2}
Assume that $f$ satisfies {\bf 3} Arguing as in the proof of
{\bf 3.} $\Rightarrow $ {\bf 1}, we obtain
the meromorphic function $F(z)$ given by (\ref{2.11}) such that
$F(z)=f(z)$ for
$z\in\D\setminus(\cZ\cup\cQ)$ (in the notation of the proof of
 {\bf 3} $\Rightarrow $ {\bf 1}). By
the already proved part of Theorem \ref{T:1.4}, we know that $f\in
\cS_{\kappa}$. By the second statement of
Lemma \ref{L:last1}, $\tilde{F}\in\cS_{\kappa'}$ for some $\kappa'\leq
\kappa +n$, where $\tilde{F}$ is the restriction of $F$ to the set
$\D\setminus(\cZ\cup\cQ)$. By continuity (statement 3. of Lemma \ref{L:new}),
the function $F$, considered on its natural domain of definition
$\D\setminus P(F)$, where $P(F)$ is the set of poles of $F$, also
belongs to $\cS_{\kappa'}$. Using Theorem \ref{T:1.2b}, write
$F=\frac{S}{B}$, where $S$ is a Schur function and $B$ is a Blaschke
product of degree $\kappa'$ without common zeros. Thus, $f$ admits an
extension to a standard function $\tilde{f}$ with $\kappa'$ poles and
$\rho$ jumps, for some nonnegative integer $\rho$. By
Theorem \ref{T:3.1}
and Remark \ref{R:last} we have
\begin{equation}\label{aga}
{\bf k}_m (f)=\kappa'+\rho, \qquad
{\rm for} \ \ {\rm every} \ \ m\geq \kappa'+ 2\rho. \end{equation}
On the other hand, since $f\in\cS_{\kappa}$, we have ${\bf k}_n(f)=\kappa$
for all sufficiently large $n$. Comparing with (\ref{aga}), we see that
$\kappa' +\rho=\kappa$, and {\bf 2} follows.
\smallskip

{\bf Proof of Theorem \ref{T:1.new}:}
If $f\in \cS_{\kappa}$, and $\tilde{f}$ is the standard function that extends
$f$ (as in statement {\bf 2} of Theorem \ref{T:1.4}), then
by Theorem \ref{T:3.1}
$$ N(f)\leq N(\tilde{f})\leq 2\kappa.$$
If
${\bf k}_{2\kappa}(f)={\bf k}_{2\kappa+3}(f)=\kappa$, then obviously {\bf 3}
of Theorem \ref{T:1.4} holds, and $f\in \cS_{\kappa}$ by Theorem \ref{T:1.4}.

\section{A local result}
\setcounter{equation}{0}

Let $\Omega\subseteq\D$ be an open set, and let
$$
{\bf k}_n(f;\Omega):=\max_{z_1,\ldots,z_n\in\Omega}\sqminus P_n(f; \,
z_1,\ldots,z_n),
$$
where the domain of definition of the function $f$ contains
$\Omega$.
A known local result
(see \cite[Theorem 1.1.4]{ADRS}, where it is proved in the more general
setting of operator valued functions)
states that for a meromorphic
function $f$ from
the class $\cS_\kappa$ and for every open set $\Omega\in\D$
that does not contain any poles of $f$, there is
an integer $n$ such that ${\bf k}_n(f;\Omega)=\kappa$. However, the minimal
such $n$
can be arbitrarily large, in contrast with Theorem \ref{T:3.1}, as the
following example shows.
\begin{Ex}
{\rm  Fix a positive integer $n$, and let
$$f_n(z)={\displaystyle \frac{z^n(2-z)}{2z-1}}
=\frac{S(z)}{b(z)}, $$
where $S(z)=z^n$ is a Schur function, and $b(z)=\frac{z-1/2}{1-z/2}$ is a
Blaschke factor. By Theorem \ref{T:1.2b},
$f_n\in \cS_1$.
Thus, the kernel
$$
K(z,w):=\frac{1-f_n(z)f_n(w)^*}{1-z{w}^*}=1+z{w}^*+\ldots
+z^{n-1}({w}^*)^{n-1}-\frac{3z^n(w^*)^n}{(2z-1)(2{w}^*-1)}
$$
has one negative square on $\D\setminus\{\frac{1}{2}\}$.
Select $n$ distinct points $\cZ=\{z_1, \ldots ,z_n\}$ (in some order) near
zero, and multiply
the Pick matrix $P_n(f;z_1, \ldots ,z_n)$
by $\Phi(\cZ)$ on the left and by $\Phi(\cZ)^*$ on the right, where
$\Phi(\cZ)$ is defined in (\ref{x.1}).
By Lemma \ref{L:x.1},
\begin{equation}\label{xyz}
\lim_{z_1, \ldots ,z_n\rightarrow 0} \Phi(\cZ)
P_n(f;z_1, \ldots ,z_n) \Phi(\cZ)^*=
\left[\frac{1}{i!}\frac{1}{j!} \left(\frac{\partial^{i+j}}{\partial z^i\partial
({w}^*)^j} K(z,w)\right)_{z,w=0}\right]_{i,j=0}^{n-1},
\end{equation}
which is the identity matrix. Therefore, by Lemma \ref{L:new}
there exists a $\delta_n>0$ such that
$P_n(f;z_1, \ldots ,z_n)\geq 0$ if $|z_1|, \ldots, |z_n|< \delta_n$.
\smallskip

On the
other hand, since
$$ \left(\frac{\partial^{2n}}{\partial z^n\partial
({w}^*)^n} K(z,w)\right)_{z,w=0}=-3, $$
selecting an ordered set of $n+1$ distinct points $\cZ'=\{z_1, \ldots
,z_{n+1}\}$ in a neighborhood
of zero, analogously to (\ref{xyz}) we obtain
\begin{eqnarray*}
\lim_{z_1, \ldots ,z_{n+1}\rightarrow 0} \Phi(\cZ')
P_{n+1}(f;z_1, \ldots ,z_{n+1}) \Phi(\cZ')^*&=&
\left[\frac{1}{i!}\frac{1}{j!} \left(\frac{\partial^{i+j}}{\partial z^i\partial
({w}^*)^j}
K(z,w)\right)_{z,w=0}\right]_{i,j=0}^{n}
\\ &
=&\left[\begin{array}{cc}I_n &0 \\ 0 &
-3\end{array}\right]. \end{eqnarray*}
By Lemma \ref{L:new} again, there exists a $\delta'_n>0$ such that
$P_{n+1}(f;z_1, \ldots ,z_{n+1})$ has exactly one negative eigenvalue if
$|z_1|, \ldots, |z_{n+1}|< \delta'_n$. (Note that
$P_{n+1}(f;z_1, \ldots ,z_{n+1})$ cannot have more than one negative eigenvalue
because $f\in \cS_{1}$.) }
\label{E:3.1a}
\end{Ex}
Theorem \ref{T:1.4}, or more precisely its proof, allows us to extend the local
result to
the more general classes
$\cS_{\kappa}$:
\begin{Tm} If $f\in \cS_{\kappa}$ is defined
on $\D\setminus \Lambda$, where $\Lambda$ is a discrete set, then for every
open set $\Omega\subseteq \D\setminus\Lambda$ there exists a positive integer
$n$
such that
\begin{equation}\label{5.z}
 {\bf k}_{n}(f;\Omega)=q + \ell_{\Omega},
\end{equation}
where
$\ell_{\Omega}$ is the
number of
jumps points of $f$
that belong to $\Omega$,
 and $q$ is the number of poles of $f$ in $\D$, counted with multiplicities.
\end{Tm}
Note that in view of {\bf 1} $\Leftrightarrow$ {\bf 2} of Theorem
\ref{T:1.4}, the right hand side of (\ref{5.z}) is equal to
$\kappa -\ell_{\not\in\Omega}$,
where $\ell_{\not\in\Omega}$
is the
number of
jump points of $f$
that do not belong to $\Omega$. Therefore,
since  ${\bf k}_{n}(f;\Omega)$ is completely
independent of the values of $f$ at jump points outside $\Omega$,
it is easy to see that
\begin{equation}\label{5.z0}
 {\bf k}_{m}(f;\Omega)\leq q + \ell_{\Omega} \qquad
{\rm for} \ \ {\rm  every}  \ \ {\rm integer}\ \  m>0. \end{equation}

{\bf Proof}: We may assume by Theorem \ref{T:1.4} that $f$ is a standard
function. We may also assume that $f$ has no jumps outside $\Omega$,
because the values of $f$ at jump points outside $\Omega$ do not
contribute to
${\bf k}_{n}(f;\Omega)$.
\smallskip

Let $w_1, \ldots ,w_k$ be the distinct poles of $f$ in $\Omega$ (if $f$ has no
poles in $\Omega$, the subsequent argument is simplified accordingly), of
orders $r_1, \ldots, r_k$, respectively, and let $z_1, \ldots, z_{\ell}$
be the jumps of $f$. Then analogously to (\ref{x.5}) we have
\begin{equation}
\label{5.z1}
f(z)=\left\{\begin{array}{cl} {\displaystyle \frac{\widehat{S}(z)}{b(z)}} &
{\rm if}\, z\not\in \{z_1,\ldots,z_{\ell}, w_{t+1},\ldots ,w_k\} \\
f_j & {\rm if }\, z=z_j, j=1, \ldots , \ell \end{array}\right.
\end{equation}
where $b(z)$
is the Blaschke product
having zeros $w_1, \ldots ,w_k$
of
orders $r_1, \ldots, r_k$, respectively.
given by (\ref{x.4b}). We assume that
$w_j=z_j$ for  $j=1, \ldots , t$, $\{w_{t+1}, \ldots ,w_k\}\cap
\{z_{t+1},\ldots
,z_{\ell}\}=\emptyset$, and $\widehat{S}(z_j)/b(z_j)\neq f_j$ for $j=t+1,
\ldots
,\ell.$ The function $\widehat{S}(z)$ is of the form
$\widehat{S}(z)=S(z)/b_{{\rm
out}}(z)$,
where $S(z)$ is a Schur function that does not vanish at $w_1, \ldots ,w_k$,
and $b_{{\rm out}}(z)$ is the Blaschke product whose zeros coincide with the
poles of $f$ outside $\Omega$, with matched orders.
 Let $N=2\ell+\sum_{j=1}^kr_j$, and select $N$ distinct points
as in (\ref{x.5a}), with the additional proviso that all these points belong to
$\Omega$.
Select also $n$ distinct points $\Xi=\{\xi_1,\ldots ,\xi_n\}$,
$\xi_j\in\Omega$,
disjoint from (\ref{x.5a}), in a vicinity of some $z_0\in \Omega$; we
assume that $f$ is analytic at $z_0$. The number $n$ of
 points $\xi_j$, and further specifications concerning the set $\Xi$, will be
determined later.
\smallskip

Let
$$ P=P_{N+n}(f; \mu_{j,i}, \nu_1, \ldots, \nu_\ell, z_1, \ldots, z_\ell, \xi_1,
\ldots , \xi_n ). $$
We now repeat the steps between (\ref{x.5a}) and (\ref{x.20}) of the proof of
Theorem \ref{T:3.1}, applied to the top left $N \times N$ corner of $P$.
As a result, we obtain the following matrix:
\begin{equation}\label{5.z2}
P_6=\left[\begin{array}{cccc}
-K_{11} & 0 & K_{31}^* & K_{41}^*\\
0& K_{22} &K_{22} & K_{42}^* \\
K_{31} & K_{22} &
K_{33} & K_{43}^* \\
K_{41}&K_{42} & K_{43} & K_{44} \end{array}\right], \end{equation}
where $K_{11}$, $K_{22}$, and $K_{31}$
are as in (\ref{x.20a}), and
\begin{eqnarray*}
K_{41}&=&\left[{\displaystyle-\frac{\widehat{S}(\xi_j)}{1-\xi_jz_{t+i}^*}}
\right]_{j,i=1}^{n,\ell-t}, \quad
K_{42}=\left[\frac{1}{1-\xi_jw^*_i}\right]_{j,i=1}^{n,\ell-t},
\end{eqnarray*}
\begin{eqnarray*}
K_{43}&=&\left[{\displaystyle\frac{1-\widehat{S}(\xi_j)\widehat{S}(w_i)^*}
{1-\xi_jw^*_i}}\right]_{j,i=1}^{n,\ell-t}, \quad
K_{44}=\left[{\displaystyle\frac{1-\widehat{S}(\xi_j)\widehat{S}(\xi_i)^*}
{1-\xi_j\xi^*_i}}\right]_{j,i=1}^{n}, \\[3mm]
K_{33}&=&\left[{\displaystyle\frac{1-
\widehat{S}(w_j)\widehat{S}(w_i)^*}{1-w_jw^*_i}}
\right]_{j,i=1}^{t}.
\end{eqnarray*}
The sizes of matrices $K_{11}$, $K_{22}$, $K_{33}$, and $K_{44}$ are
$(\ell -t)\times (\ell -t)$, $t \times t$, $t \times t$, and $n \times n$,
respectively.
It follows that
\begin{equation}\label{5.z2a}
\sqminus P \geq N-2\ell+ \sqminus P_6= \sum_{j=1}^k r_j + \sqminus P_6.
\end{equation}
Take the Schur complement ${\bf S}$
to the block
$\left[\begin{array}{cc} -K_{11}&0\\ 0 &K_{22} \end{array}\right]$ in
(\ref{5.z2}):
$$ {\bf S}=\left[\begin{array}{cc} {\bf S}_{11} &
{\bf S}_{12} \\ {\bf S}_{21} & {\bf S}_{22}\end{array}\right],
$$
where
$$
{\bf S}_{22}=K_{44}-K_{42}K_{22}^{-1}K_{42}^*+K_{41}K_{11}^{-1}K_{41}^*.
$$
We have
\begin{equation}\label{5.z2b}
\sqminus P_6= \ell -t +\sqminus {\bf S}. \end{equation}
Analogously to formulas (\ref{x.22}) and (\ref{x.23}) we obtain
\begin{equation}
\label{5.z4}
{\bf S}_{11}=
-\left[\frac{\widehat{S}(w_j)\vartheta(w_j)\vartheta(w_i)^*
\widehat{S}(w_i)^*}{1 -w_jw_i^*}\right]_{j,i=1}^t,\quad
{\bf S}_{21}=
-\left[\frac{\widehat{S}(\xi_j)\vartheta(\xi_j)\vartheta(w_i)^*
\widehat{S}(w_i)^*}{1-\xi_jw_i^*}\right]_{j,i=1}^{t,n},
\end{equation}
 where the rational function $\vartheta$ is given by (\ref{x.20b}).
Note that $\widehat{S}(w_i)\vartheta(w_i)\neq 0$, for $i=1, \ldots ,t$.
Multiply ${\bf S}$ by the matrix
\begin{equation}
\label{5.z5}
{\rm diag}\, \left(\frac{1}{\widehat{S}(w_1)\vartheta(w_1)}, \ldots ,
\frac{1}{\widehat{S}(w_t)\vartheta(w_t)}, I\right) \end{equation}
on the left and by the adjoint of (\ref{5.z5}) on the right. We arrive at
\begin{equation}\label{5.z6}
P_7=\left[\begin{array}{cc}
-\left[{\displaystyle\frac{1}{1-w_jw_i^*}}\right]_{j,i=1}^t &
\left(-\left[{\displaystyle\frac{\widehat{S}(\xi_j)\vartheta(\xi_j)}{1-
\xi_jw_i^*}}\right]_{j,i=1}^{n,t}\right)^*\\
-\left[{\displaystyle\frac{\widehat{S}(\xi_j)\vartheta(\xi_j)}{1-
\xi_jw_i^*}}\right]_{j,i=1}^{n,t} &
K_{44}-K_{42}K_{22}^{-1}K_{42}^*+K_{41}K_{11}^{-1}K_{41}^* \end{array}\right].
\end{equation}
Finally, we take the Schur complement, denoted $P_8$, to the block
$-\left[{\displaystyle \frac{1}{1-w_jw_i^*}}\right]_{j,i=1}^t$ of $P_7$.
Then
\begin{equation}\label{5.z2c}
\sqminus {\bf S}=t+ \sqminus P_8. \end{equation}
The $(j,i)$ entry
of $P_8$ is
\begin{eqnarray*}
&&\frac{1-\widehat{S}(\xi_j)\widehat{S}(\xi_i)^*}{1-\xi_j\xi_i^*} -G_{\ell
-t}^*(I-\xi_jD(\cW_1)^*)^{-1}K_{22}^{-1}(I-\xi_i^*D(\cW_1))^{-1}G_{\ell-t}\\
&&+ \widehat{S}(\xi_j)G_{\ell-t}^*
(I-\xi_jD(\cZ_2)^*)^{-1}K_{11}^{-1}(I-\xi_i^*D(\cZ_2))^{-1}G_{\ell-t}
\widehat{S}(\xi_i)^* \\
&&+\widehat{S}(\xi_j)\vartheta(\xi_j)
G_{\ell  -t}^*(I-\xi_jD(\cW_1)^*)^{-1}K_{22}^{-1}
(I-\xi_i^*D(\cW_1))^{-1}G_{\ell-t}\vartheta(\xi_i)^*\widehat{S}(\xi_i)^*.
\end{eqnarray*}
Since
$$
K_{22} -D(\cW_1)K_{22}D(\cW_1)^*=G_{\ell-t}G_{\ell -t}^*,
$$
we conclude (as it was done for the function $\vartheta$ given in
(\ref{x.20b})) that the rational function
$$
\widehat{\vartheta}(z)=1 +(z-1) G^*_{\ell -t}
(I-zD(\cW_1)^*)^{-1}K_{22}^{-1}
(I-D(\cW_1))^{-1}G_{\ell-t}
$$
is inner, and its set of zeros is $Z(\widehat{\vartheta})=\cW_1$.
Using (\ref{x.22a}) and  the formula
$$
1-\widehat{\vartheta}(z)\widehat{\vartheta}(w)^*=G^*_{\ell -t}
(I-zD(\cW_1)^*)^{-1}K_{22}^{-1}(I-w^*D(\cW_1))^{-1}G_{\ell-t}
$$
similar to (\ref{x.22a}), the expression for the $(j,i)$ entry of $P_8$ takes
the form
\begin{eqnarray*}
&&\frac{1-\widehat{S}(\xi_j)\widehat{S}(\xi_i)^*}{1-\xi_j\xi_i^*}-
\frac{1-\widehat{\vartheta}(\xi_j)\widehat{\vartheta}(\xi_i)^*}
{1-\xi_j\xi_i^*}+\widehat{S}(\xi_j)
\frac{1-\vartheta(\xi_j)\vartheta(\xi_i)^*}{1-\xi_j\xi_i^*}
\widehat{S}(\xi_i)^*\\
&&+\widehat{S}(\xi_j)\vartheta(\xi_j)
\frac{1-\widehat{\vartheta}(\xi_j)\widehat{\vartheta}(\xi_i)^*}{1-\xi_j\xi_i^*}
\vartheta(\xi_i)^*\widehat{S}(\xi_i)^*\\
&&= \widehat{\vartheta}(\xi_j)
\frac{1-
\widehat{S}(\xi_j)\vartheta(\xi_j)\vartheta(\xi_i)^*\widehat{S}(\xi_i)^*}
{1-\xi_j\xi_i^*}\widehat{\vartheta}(\xi_i)^*.
\end{eqnarray*}
We now assume that the point $z_0\in\Omega$ in a neighborhood of which the set
$\Xi$ is selected is such that $\widehat{\vartheta}(z_0)\neq 0$.
Then we may assume that $\widehat{\vartheta}(\xi_j)\neq 0$, $j=1, \ldots,
n$. Multiply $P_8$ on the left by the matrix
${\rm diag}\, (\widehat{\vartheta}(\xi_1)^{-1}, \ldots ,
\widehat{\vartheta}(\xi_n)^{-1})$ and on the right by
its adjoint of, resulting in a matrix
$$
P_9=\left[\frac{1-T(\xi_j)T(\xi_i)^*}{1-\xi_j\xi_i^*}\right]_{j,i=1}^n,
$$
where the function $T(z)$ is given by $T(z)=\widehat{S}(z)\vartheta(z)$.
We have
\begin{equation}\label{5.z2d}
\sqminus P_8=\sqminus P_9 \end{equation}
Note that $T(z)$ is a meromorphic function with no poles in $\Omega$.
By \cite[Theorem 1.1.4]{ADRS}, there exist a positive integer $n$ and
points $\xi_1, \ldots , \xi_n$ in a neighborhood of $z_0$ such that
$P_9$ has $q - \sum_{j=1}^k r_j$ (the number of poles of $T(z)$) negative
eigenvalues.
Using this choice of $\xi_j$, and
combining (\ref{5.z2a}), (\ref{5.z2b}), (\ref{5.z2c}), and
(\ref{5.z2d}), we obtain
$$ \sqminus P\geq q+\ell. $$
Now equality (\ref{5.z}) follows in view of (\ref{5.z0}).
\qed

\section{An open problem}
\setcounter{equation}{0}

Fix an integer $k\geq 3$.
An open set $D\subseteq \D$ is said to have the $k$-{\em th
Hindmarsh property} if every function $f$ defined on $D$ and such that
the matrices $P_k(f;z_1,\ldots ,z_k)$
are positive semidefinite for every $k$-tuple of points
$z_1, \ldots ,z_k\in D$, admits a (necessarily unique) extension
to a Schur function (defined on $\D$). Theorem \ref{T:1.2}
shows that
$\D$ has the $3$-rd Hindmarsh property.
Example \ref{E:3.1a} shows that the open disk of radius $\delta_n$ centered at
the origin
does not have the $n$-th Hindmarsh property, if $\delta_n$ is sufficiently
small.
\medskip

Denote by ${\cH}_k$ the collection of all sets with the $k$-th
Hindmarsh
property. Clearly,
$$ {\cH}_3 \subseteq {\cH}_4 \subseteq \cdots \subseteq {\cH}_k
\subseteq \cdots  $$
\begin{Pn}
Let $D\in {\cH}_k$. If $\cZ\subset D$ is a discrete set in $D$, then
$D\setminus \cZ$ also has the $k$-th Hindmarsh property.
\end{Pn}
{\bf Proof.} Let $f$ be a function with the domain of definition
$D\setminus \cZ$ and such that $P_k(f;z_1,\ldots ,z_k)\geq 0$
for every set of distinct points $z_1, \ldots ,z_k\in D\setminus \cZ$.
In particular,
$P_3(f;z_1,z_2,z_3)\geq 0$
for every triple of distinct points $z_1, z_2,z_3\in D\setminus \cZ$.
By Hindmarsh's theorem, $f$ is analytic on $D\setminus \cZ$.
Also, $$|f(z)|^2=(|f(z)|^2-1)+1=-P_1(f;z)(1-|z|^2)+1 \leq 1 $$
for every $z\in D\setminus \cZ$. Thus, $f$ admits an analytic
continuation to a function $\widehat{f}$ on $D$. By continuity,
$P_k(\widehat{f};z_1,\ldots ,z_k)\geq 0$ for every $k$-tuple of
distinct points $z_1,\ldots ,z_k \in D$. Since $D\in {\cH}_k$,
$\widehat{f}$ admits an extension to a Schur function.
\qed
\bigskip

\begin{Cy} Let $D_0=\D$, $D_j=D_{j-1}\setminus \cZ_{j-1}$, $j=1,2,
\ldots$, where $\cZ_{j-1}$ is a discrete set in $D_{j-1}$. Then all the
sets $D_j$, $j=1,2, \ldots, $, have the $3$-rd Hindmarsh property.
\end{Cy}

Using the sets with the Hindmarsh property,
the implication {\bf  3.} $\Rightarrow$ {\bf 1.} of Theorem \ref{T:1.4}
can be extended to a larger class of functions, as follows:
\begin{Tm}
Let $f$ be defined on $D$, where $D\in {\cH}_p$. Then
$f$ belongs to $\cS_\kappa (D)$ if and only if
\begin{equation}
{\bf k}_n(f)={\bf k}_{n+p}(f)=\kappa
\label{1.8'}
\end{equation}
for some integer $n$. \label{T:last}
\end{Tm}
The definition of $\cS_\kappa (D)$ is the same as $\cS_\kappa$, with the only
difference being that $\cS_\kappa (D)$ consists of functions defined on $D$.
The proof of Theorem \ref{T:last} is essentially the same as that of
{\bf 3.} $\Rightarrow$ {\bf 1.} of Theorem \ref{T:1.4}.
\begin{Pb}
Describe the structure of sets with the $k$-th Hindmarsh property.
\end{Pb}
In the general case this seems to be a very difficult unsolved problem.
We do not address this problem in the present paper.
\bigskip

{\bf Acknowledgement}. The research of LR was supported in
part by an
NSF grant.
\bigskip

\bibliographystyle{amsplain}
\providecommand{\bysame}{\leavevmode\hbox to3em{\hrulefill}\thinspace}

\parbox[t]{7 cm}{
Department of Mathematics \\
P. O. Box 8795 \\
The College of William and Mary\\
Williamsburg VA 23187-8795, USA \\
vladi@@math.wm.edu, sykhei@@wm.edu, lxrodm@@math.wm.edu }

\end{document}